\documentclass{article}%
\usepackage{amsmath}
\usepackage{amsfonts}
\usepackage{amssymb}
\usepackage{graphicx}%
\providecommand{\U}[1]{\protect\rule{.1in}{.1in}}
\newtheorem{theorem}{Theorem}

\newtheorem{corollary}{Corollary}

\newtheorem{lemma}{Lemma}

\begin{document}

\title{Optimal robust estimators for families of distributions on the
integers\thanks{\noindent\textbf{Acknowledgement.} This research was partially
supported by Grants X-094 and 20020170100022BA from Universidad de Buenos
Aires, PID 5505 from CONICET and PAV 120 and PICT \ 21407 from ANPCYT,
Argentina. {}}}
\author{Ricardo A. Maronna$^{1}$\thanks{Corresponding author. E-mail
rmaronna@retina.ar} and Victor J. Yohai$^{2}$\\\ $^{1}$Universidad Nacional de La Plata\\$^{2}$Universidad de Buenos Aires and CONICET}
\date{}
\maketitle

\begin{abstract}
Let $F_{\theta}$ be a family of distributions with support on the set of
nonnegative integers $Z_{0}$. In this paper we derive the M-estimators with
smallest gross error sensitivity (GES). We start by defining \ \ the uniform
median \ of a distribution $F$ with support on $Z_{0}$ (umed$(F)$) as the
median of $x+u,$ where $x$ and $u$ \ are independent variables with
distributions $F$ \ and uniform in [-0.5,0.5] respectively.\ Under some
general conditions we prove that the estimator with smallest GES satisfies
umed$(F_{n})=$umed$(F_{\theta}),$ where $F_{n}$ is the empirical distribution.
\ The asymptotic distribution of these estimators is found. This distribution
is normal \ except when there is a positive integer $k$ \ so that $F_{\theta
}(k)=0.5.$ In this last case, the asymptotic distribution behaves as normal at
each side of 0, but with different variances. A simulation Monte Carlo study
compares, for the Poisson distribution, the efficiency and robustness for
finite sample sizes of this estimator with those of other robust estimators.


\end{abstract}

\section{Introduction}

Consider a one-parameter family of distributions $F_{\theta}.$ An important
problem in the theory of robust estimation is the study of estimators which in
some sense optimize their bias under contamination. The \emph{gross-error
sensitivity} (GES) is defined as the maximum of the absolute values of the
influence function. It gives an approximation to the maximum bias produced by
an outlier contamination of rate $\varepsilon,$ when $\varepsilon$ is
\textquotedblright small\textquotedblright.

Hampel (1974) dealt with M-estimators defined as solutions of equations of the
form%
\[
\sum_{i=1}^{n}\psi\left(  x_{i},\theta\right)  =0,
\]
and considered the problem of minimizing the asymptotic variance among
Fisher-consistent M-estimators which satisfy a bound on the GES.
Alternatively, this problem can be stated as minimizing the GES under a bound
on the asymptotic variance. Details are given in Section \ref{secHampel}.

In this paper we consider minimizing the GES without any restrictions on the
asymptotic variance. Let $F_{\theta}$ be a family of continuous distributions
with densities $p(x,\theta)$ and score function $\ $%
\begin{equation}
\psi_{0}(x,\theta)=\frac{\partial\log p(x,\theta)}{\partial\theta}.
\label{defPsi0}%
\end{equation}

Maronna et al. (2019, p. 68) show that if $\psi_{0}(x,\theta)$ is strictly
monotone on $x,$ the M-estimator with smallest GES is obtained \ by solving
\begin{equation}
\text{\textrm{med}}\left(  F_{n}\right)  \text{=\textrm{med}}\left(
F_{\widehat{\theta}}\right)  \text{,} \label{medmed}%
\end{equation}
where \textquotedblleft med\textquotedblright\ stands for the median and
$F_{n}$ denotes the empirical distribution.

This result does not hold when $F_{\theta}$ has support on a discrete set,
such as the integers, since\ \textrm{med}$(F_{\theta})$ is in general not
uniquely defined and therefore (\ref{medmed}) does not identify $\theta.$

To overcome this problem, we introduce \ in Section \ref{secUmed} \ the
\ concept of the \emph{uniform median} of a distribution $F$ with support on
the set of nonnegative integers $Z_{0}.$ In Section \ref{secOptimal} we prove
that the estimator with smallest GES \ can be obtained by solving \ an
equation similar to (\ref{medmed}) but replacing the median by the uniform
median. Section \ref{secAsympto} deals with the asymptotic distribution of
this estimator. Section \ref{secPoisson} shows an application to the family of
Poisson distributions. Finally Section \ref{SecApend} is an appendix
containing proofs of the main results.

\section{The uniform median\label{secUmed}}

We shall deal with distributions concentrated on a finite or infinite interval
of the integers, such as the Poisson distribution. To avoid notational
complications it will henceforth be assumed that this interval is the set
$Z_{0}$ of the nonnegative integers. Let $\ F$ be a distribution with support
on $Z_{0}$ and call $p$ the corresponding probability density. The
\emph{uniform median} of $F$ (\textrm{umed}$(F)$) is defined as the median of
the distribution that\ distributes the mass $p(k)$ uniformly on the interval
$[k-0.5,k+0.5].$ This is equivalent\ to define \textrm{umed}$(F)$ as the
median of $x+u$ where $x$ and $\ u$ are independent with distributions $F$ and
uniform in $[-0.5,0.5]$ respectively. To give an explicit formula for
\textrm{umed}$(F),$ define%
\begin{equation}
k_{0}(F)=\min\{k:F(k)\geq0.5\}. \label{defK0}%
\end{equation}
It is easy to verify that
\begin{equation}
\text{\textrm{umed}}(F)=k_{0}(F)-0.5+\frac{0.5-F(k_{0}(F)-1)}{p_{0}},
\label{umeddef}%
\end{equation}
where%
\begin{equation}
p_{0}=p(k_{0}(F)). \label{def-p0}%
\end{equation}
Note that by definition $p_{0}>0,$ and therefore \textrm{umed}$(F)$ is well defined.

Ma et al. (2011) define quantiles of discrete distributions based on what they
call \textquotedblleft mid-distribution functions\textquotedblright. In the
case of the median, their definition is similar to that of the umedian but is
not exactly equal, nor can any of the two medians be expressed as a function
of the other one.

The following property is immediate%

\[
k_{0}(F)-0.5<\text{\textrm{umed}}(F)\leq k_{0}(F)+0.5.
\]
Therefore for two distributions $F_{1}$ and $F_{2}$ we have%
\begin{equation}
\text{\textrm{umed}}(F_{1})=\text{\textrm{umed}}(F_{2})\Rightarrow k_{0}%
(F_{1})=k_{0}(F_{2}). \label{impk0}%
\end{equation}

Recall that the median of $F$ is the solution $\mu$ of
\[
\mathrm{E}_{F}\text{\textrm{sign}}(x-\mu)=0\text{.}%
\]
Similarly \textrm{umed}$(F)$ can be defined as the solution $\mu$ of
\[
\mathrm{E}\psi_{0.5}^{H}(x-\mu)=0,
\]
where $\psi_{m}^{H}$ is the Huber family of score functions given by%
\begin{equation}
\psi_{m}^{H}(x)=\max(\min(x,m),-m), \label{defHuber}%
\end{equation}

Note that $\psi_{0.5}^{H}(x)=0.5$sign$(x)$ for $|x|\geq0.5$ and $\psi
_{0.5}^{H}(x)=x$ for $|x|<0.5.$

\section{Estimators with smallest GES\label{secOptimal}}

\subsection{The Hampel approach\label{secHampel}}

Before dealing with the unrestricted minimization of the GES we need to
consider Hampel's (1974) approach to robust optimality.\ These optimal
estimators are defined by minimizing the asymptotic variance among the class
of M-estimators which are Fisher-consistent and have their GES bounded by a
given constant. \ 

If $x_{1},...,x_{n}$ is a a random sample from $F_{\theta}$, the
Hampel-optimal estimator $\widehat{\theta}$ is defined as the solution of the
estimating equation%
\begin{equation}%
{\displaystyle\sum_{i=1}^{n}}
\psi_{m}^{H}[\psi_{0}(x_{i},\widehat{\theta})-c(m,\widehat{\theta})]=0,
\label{esteq}%
\end{equation}
where $m$ \ depends on the given bound on the GES, $\psi_{0}$ and $\psi
_{m}^{H}$ are defined in (\ref{defPsi0}) and (\ref{defHuber}), respectively,
and $c(m,\theta)$ is defined \ by
\begin{equation}
\mathrm{E}_{\theta}\psi_{m}^{H}(\psi_{0}(x_{i},\theta)-c(m,\theta))=0,
\label{conHamop}%
\end{equation}
where $F_{\theta}$ denotes the expectation with respect to $F_{\theta}.$

The \emph{dual} Hampel problem consists of minimizing the GES under a bound
$V$ on the asymptotic variance. It is known that the solution has again the
form (\ref{esteq}), where $m$ is a decreasing function of $V.$

\subsection{The optimal estimator}

Let $F_{\theta}$ be a family of distribution functions with $\theta\in
\Theta\subset R,$ with support on $Z_{0}$ and probability densities
$p(x,\theta)$. \ Given \ a random sample $x_{1},...x_{n}$ of $F_{\theta},$
denote by $F_{n}$ the corresponding empirical distribution function, with
density \ \
\[
p_{n}(k)=\frac{\#\{i:x_{i}=k\}}{n}.
\]
Then we have:

\begin{theorem}
\label{Theo-Opti}\textbf{\ }Assume that $\psi_{0}(x,\theta)$ is continuous and
strictly monotone in $x$ and $\theta.$ Then the estimator with smallest GES is
$\widehat{\theta}_{n}$ defined by%
\begin{equation}
\text{\textrm{umed}}(F_{n})=\text{\textrm{umed}}(F_{\widehat{\theta}_{n}}).
\label{eqmed}%
\end{equation}

\end{theorem}

The proof is given in Section \ref{Apend-Opti}

Before showing the existence and uniqueness of the solution to (\ref{eqmed})
we state some properties of the uniform median.

It will henceforth be assumed that

\begin{description}
\item[A1] $F_{\theta}\left(  x\right)  $ is a decreasing function of $\theta$
for all $x$

\item[A2] $F_{\theta}\left(  x\right)  $ is a continuous function of $\theta$
for all $x$

\item[A3] Call $\Theta$ the range of $\theta$ (e.g., $[0,\infty)$ for the
Poisson distribution) and let $\theta_{1}=\inf\left(  \Theta\right)  $ and
$\theta_{2}=\sup\left(  \Theta\right)  .$ Then $\lim_{\theta\rightarrow
\theta_{1}}F_{\theta}\left(  x\right)  =0$ for $x\neq0$ and $\lim
_{\theta\rightarrow\theta_{2}}F_{\theta}\left(  x\right)  =0$ for all $x$
\end{description}

These assumptions are satisfied by the standard discrete families such as the
Poisson family. A1 implies that the family is \textquotedblleft stochastically
increasing\textquotedblright\ in the sense that $\theta_{1}>\theta_{2}$
implies that $F_{\theta_{1}}$ is strictly stochastically larger than
$F_{\theta_{2}}$.

\begin{lemma}
\label{Lema-Lim}Let $F^{(n)},$ $n\geq1$ be a sequence of distributions with
support on $Z_{0}$ . Then $F^{(n)}\rightarrow_{w}F$ (where $\rightarrow_{w}%
$denotes weak convergence) implies that \textrm{umed}$(F^{\left(  n\right)
})\rightarrow$ \textrm{umed}$(F).$
\end{lemma}

The proof is given in Section \ref{Apend-LemaLim}.

\begin{corollary}
Let $x_{1},...,x_{n}$ be i.i.d random variables with distribution $F$ with
support \ at $Z_{0}$ and call\ $F_{n}$ the empirical distribution.\ Then
\textrm{umed}$(F_{n})$ $\rightarrow$\textrm{umed}$(F)$ a.s.$.$
\end{corollary}

\textbf{Proof:\ }The result follows from Lemma \ref{Lema-Lim} and the fact
that $F_{n}\rightarrow_{w}F$ by the Glivenko-Cantelli Theorem.

\begin{lemma}
Put for brevity $g\left(  \theta\right)  =$\textrm{umed}$\left(  F_{\theta
}\right)  .$ If A1-A2-A3 hold, then $g\left(  \theta\right)  $ is a continuous
increasing function of $\theta$, and $\lim_{\theta\rightarrow\theta_{1}%
}g\left(  \theta\right)  =0$ and $\lim_{\theta\rightarrow\theta_{2}}g\left(
\theta\right)  =\infty.$
\end{lemma}

\textbf{Proof: } A1 implies that $g$ is increasing; its continuity follows
from A2 and Lemma \ref{Lema-Lim}, and A3 implies the last statement.

\begin{corollary}
(\ref{eqmed}) has a unique solution.
\end{corollary}

\section{Asymptotics\label{secAsympto}}

\subsection{Consistency.}

The following result proves the strong consistency of the optimal estimator
$\widehat{\theta}_{n}$ defined by (\ref{eqmed}).

\begin{theorem}
\label{cons}\textbf{\ }Let $F_{\theta}$ satisfy A1-A2-A3 and let
$x_{1},...,x_{n}$ be i.i.d. random variables with distribution $F_{\theta}.$
Then $\widehat{\theta}_{n}\rightarrow\theta$ a.s.
\end{theorem}

\textbf{Proof.} Call $F_{n}$ the empirical distribution. and put $g\left(
\theta\right)  =$\textrm{umed}$\left(  F_{\theta}\right)  .$ Since
$\widehat{\theta}_{n}=g^{-1}($\textrm{umed}($F_{n})),$ and $g^{-1}$ is
continuous, we have%
\begin{align*}
\lim_{n\rightarrow\infty}\widehat{\theta}_{n}  &  =g^{-1}(\lim_{n\rightarrow
\infty}\text{\textrm{umed}}(F_{n}))=g^{-1}(\lim_{n\rightarrow\infty
}\text{\textrm{umed}}(F_{\theta})).\\
&  =g^{-1}(g(\theta))=\theta\text{ a.s..}%
\end{align*}
This proves the Theorem.

\subsection{ Asymptotic distribution}

The following Theorem states that when $F(k_{0}(F))>0.5,$ \textrm{umed}%
$(F_{n})$ is asymptotically normal, while when \ $F(k_{0}(F))=0.5,$ its left
and right tails are asymptotically normal but with different variances.

\begin{theorem}
\label{Theo-Dist-Umed}$\ $ Let $x_{1},...,x_{n}$ \ be i.i.d. random variables
with distribution $F$ with support \ at $Z_{0}$, $p$ its probability density
and \ $F_{n}$ the empirical distribution.\ Put for brevity
\[
K=k_{0}\left(  F\right)  ,\ p_{0}=p\left(  K\right)  ,\ F^{1}=F\left(
K-1\right)  \ \mathrm{and}\ Z_{n}=n^{1/2}(\text{\textrm{umed}}(F_{n}%
)-\text{\textrm{umed}}(F)).
\]
Then,
\end{theorem}

\textbf{(a)} If $F(K)>0.5,$ then%
\begin{equation}
Z_{n}\rightarrow_{d}N(0,\sigma^{2}) \label{Normal(a)}%
\end{equation}
where $\rightarrow_{d}$ stands for convergence in distribution and
\begin{equation}
\sigma^{2}=\frac{0.25}{p_{0}^{3}}\left(  4F^{1}\left(  F^{1}-1+p_{0}\right)
-p_{0}+1\right)  \allowbreak. \label{sigma-umed}%
\end{equation}

(b) If $F(K)=0.5,$ then
\[
Z_{n}\rightarrow_{d}H
\]
with $\ $%
\begin{equation}
H\left(  t\right)  =\left\{
\begin{array}
[c]{ccc}%
\Phi(2tp\left(  K\right)  ) & \text{if} & a\leq0\\
\Phi(2tp\left(  K+1\right)  ) & \text{if} & a>0
\end{array}
\right.  , \label{defH}%
\end{equation}
where $\Phi$ is the standard normal distribution function.

The proof is given in Section \ref{Apend-Dist-Umed}.

\textbf{Remark:} Note that the phenomenon of an asymptotic distribution with
two different normal tails also occurs with the ordinary median of samples
from a continuous distribution if the density has different side derivatives
at the population median.

The following Theorem deals with the asymptotic distribution of the
GES-optimal estimator.

\begin{theorem}
\label{Teo_Dist-Opti}\textbf{\ }Let $x_{1},...,x_{n}$ be a \ random sample of
$F_{\theta}.$ \ Assume that the density $p(x,\theta)$ is continuously
differentiable in $\theta$. Let $\widehat{\theta}_{n}$ be the estimator
defined by (\ref{eqmed}). Let
\begin{equation}
g(\theta)=k_{0}(F_{\theta})-0.5+\frac{0.5-F_{\theta}(k_{0}(F_{\theta}%
)-1)}{p(k_{0}(F_{\theta}),\theta)}. \label{defg(theta)}%
\end{equation}

\end{theorem}

(a) If $F_{\theta}(k_{0}(F_{\theta}))>0.5$, then
\[
n^{1/2}(\widehat{\theta}_{n}-\theta)\rightarrow_{d}N\left(  0,\frac{\sigma
^{2}}{g^{\prime}(\theta)^{2}}\right)
\]
with $\sigma$ defined in (\ref{sigma-umed}). \ \ 

(b ) If $F(k_{0}(F))=0.5$ then. \
\[
n^{1/2}(\widehat{\theta}_{n}-\theta)\rightarrow_{d}G,
\]
where
\[
G\left(  t\right)  =\left\{
\begin{array}
[c]{ccc}%
\Phi(2g_{-}^{\prime}tp\left(  K,\theta\right)  ) & \text{if} & a\leq0\\
\Phi(2g_{+}^{\prime}tp\left(  K+1,\theta\right)  ) & \text{if} & a>0
\end{array}
\right.  ,
\]

with
\begin{equation}
g_{+}^{\prime}=\frac{\partial}{\partial t}\left[  \frac{0.5-F_{t}\left(
K\right)  }{p\left(  K+1,\theta\right)  }\right]  _{t=\theta+}%
\ \ \mathrm{and\ \ }g_{-}^{\prime}=\frac{\partial}{\partial t}\left[
\frac{0.5-F_{t}\left(  K-1\right)  }{p\left(  K,\theta\right)  }\right]
_{t=\theta-}. \label{laterales}%
\end{equation}

The proof is given in Section \ref{Apend-Dist-Opti}.

\section{\medskip Application to the Poisson distribution\label{secPoisson}}

In this section we compute the maximum asymptotic bias of the GES-optimal
estimator and compare it to those of two robust estimators: the MT estimator
of (Valdora and Yohai 1974) and the Quasi-Likelihood (Q-L) estimator of
(Cantoni and Ronchetti 2001). Note that the GES is only a rough measure of the
bias, and therefore minimizing the GES does not ensure any optimal properties
of the actual bias.

For the Poisson distribution with parameter $\lambda$ we contaminate the data
with a point mass with probability $\varepsilon$ located at $x_{0}.$ For each
estimator $\widehat{\lambda}$ the absolute bias $|\widehat{\lambda}-\lambda|$
is computed, and the maximum over all $x_{0}$ is reported in Table
\ref{TabBias}.

\begin{center}%
\begin{table}[tbp] \centering
\caption{Maximum asymptotic biases for the Poisson family}%
\begin{tabular}
[c]{ccccc}\hline
$\varepsilon$ & $\lambda$ & Optimal & MT & Q-L\\\hline
\multicolumn{1}{r}{0.1} & \multicolumn{1}{r}{5} & \multicolumn{1}{r}{0.329} &
\multicolumn{1}{r}{0.421} & \multicolumn{1}{r}{0.409}\\
\multicolumn{1}{r}{} & \multicolumn{1}{r}{10} & \multicolumn{1}{r}{0.511} &
\multicolumn{1}{r}{0.608} & \multicolumn{1}{r}{0.627}\\
\multicolumn{1}{r}{} & \multicolumn{1}{r}{20} & \multicolumn{1}{r}{0.823} &
\multicolumn{1}{r}{0.985} & \multicolumn{1}{r}{0.959}\\
\multicolumn{1}{r}{0.2} & \multicolumn{1}{r}{5} & \multicolumn{1}{r}{0.805} &
\multicolumn{1}{r}{1.087} & \multicolumn{1}{r}{1.071}\\
\multicolumn{1}{r}{} & \multicolumn{1}{r}{10} & \multicolumn{1}{r}{1.052} &
\multicolumn{1}{r}{1.413} & \multicolumn{1}{r}{1.402}\\
\multicolumn{1}{r}{} & \multicolumn{1}{r}{20} & \multicolumn{1}{r}{1.569} &
\multicolumn{1}{r}{2.057} & \multicolumn{1}{r}{2.071}\\\hline
\end{tabular}
\label{TabBias}%
\end{table}%

\end{center}

Table \ref{TabEfi} gives the asymptotic efficiencies of the three estimators.

\begin{center}%
\begin{table}[tbp] \centering
\caption{Asymptotic fficiencies of estimators}%
\begin{tabular}
[c]{cccc}\hline
$\lambda$ & Optimal & MT & Q-L\\\hline
\multicolumn{1}{r}{5} & \multicolumn{1}{r}{0.72} & \multicolumn{1}{r}{0.89} &
\multicolumn{1}{r}{0.93}\\
\multicolumn{1}{r}{10} & \multicolumn{1}{r}{0.69} & \multicolumn{1}{r}{0.93} &
\multicolumn{1}{r}{0.96}\\
\multicolumn{1}{r}{20} & \multicolumn{1}{r}{0.67} & \multicolumn{1}{r}{0.93} &
\multicolumn{1}{r}{0.97}\\\hline
\end{tabular}
\label{TabEfi}%
\end{table}%

\end{center}

These results show that the optimal estimator has a comparatively good bias
behavior not only for \textquotedblleft small\textquotedblright\ $\varepsilon
.$ At the same time, the price for such a low bias is a relatively low efficiency.

To study the small sample behavior of the estimator a simulation with $n=20$
and 50 was run, with 500 replications. Table \ref{TabEfiFini} shows the
estimator's efficiencies, which are seen to differ little from the asymptotic ones.

\begin{center}%
\begin{table}[tbp] \centering
\caption{Finite sample efficiencies of estimators}%
\begin{tabular}
[c]{ccccc}\hline
$n$ & $\lambda$ & Optimal & MT & Q-L\\\hline
\multicolumn{1}{r}{20} & \multicolumn{1}{r}{5} & \multicolumn{1}{r}{0.80} &
\multicolumn{1}{r}{0.95} & \multicolumn{1}{r}{0.98}\\
\multicolumn{1}{r}{} & \multicolumn{1}{r}{10} & \multicolumn{1}{r}{0.64} &
\multicolumn{1}{r}{0.87} & \multicolumn{1}{r}{0.93}\\
\multicolumn{1}{r}{} & \multicolumn{1}{r}{20} & \multicolumn{1}{r}{0.56} &
\multicolumn{1}{r}{0.83} & \multicolumn{1}{r}{0.89}\\\hline
\multicolumn{1}{r}{50} & \multicolumn{1}{r}{5} & \multicolumn{1}{r}{0.72} &
\multicolumn{1}{r}{0.89} & \multicolumn{1}{r}{0.94}\\
\multicolumn{1}{r}{} & \multicolumn{1}{r}{10} & \multicolumn{1}{r}{0.71} &
\multicolumn{1}{r}{0.92} & \multicolumn{1}{r}{0.96}\\
\multicolumn{1}{r}{} & \multicolumn{1}{r}{20} & \multicolumn{1}{r}{0.66} &
\multicolumn{1}{r}{0.91} & \multicolumn{1}{r}{0.94}\\\hline
\end{tabular}
\label{TabEfiFini}%
\end{table}%

\end{center}

To study the estimators' robustness, in each sample a proportion $\varepsilon$
of the values was replaced by a value $x_{0},$ and the estimators' means
squared error (MSE) was computed. This was done for $x_{0}=0,1,...,3\lambda.$
Table \ref{TabMSE} displays the maximum MSE of each estimator over all $x_{0}%
$'s, for $\varepsilon=0.1$ and 0.2.

\begin{center}%
\begin{table}[tbp] \centering
\caption{Maximum MSEs of estimators for contamination $\epsilon$}%
\begin{tabular}
[c]{cccccc}\hline
$n$ & $\varepsilon$ & $\lambda$ & Optimal & MT & Q-L\\\hline
20 & \multicolumn{1}{r}{0.1} & \multicolumn{1}{r}{5} &
\multicolumn{1}{r}{0.67} & \multicolumn{1}{r}{0.54} & \multicolumn{1}{r}{0.52}%
\\
& \multicolumn{1}{r}{} & \multicolumn{1}{r}{10} & \multicolumn{1}{r}{1.14} &
\multicolumn{1}{r}{1.10} & \multicolumn{1}{r}{1.04}\\
& \multicolumn{1}{r}{} & \multicolumn{1}{r}{20} & \multicolumn{1}{r}{2.53} &
\multicolumn{1}{r}{2.22} & \multicolumn{1}{r}{2.32}\\\cline{2-6}
& \multicolumn{1}{r}{0.2} & \multicolumn{1}{r}{5} & \multicolumn{1}{r}{1.22} &
\multicolumn{1}{r}{1.72} & \multicolumn{1}{r}{1.61}\\
& \multicolumn{1}{r}{} & \multicolumn{1}{r}{10} & \multicolumn{1}{r}{2.26} &
\multicolumn{1}{r}{2.81} & \multicolumn{1}{r}{2.82}\\
& \multicolumn{1}{r}{} & \multicolumn{1}{r}{20} & \multicolumn{1}{r}{4.63} &
\multicolumn{1}{r}{5.81} & \multicolumn{1}{r}{5.85}\\\hline
50 & 0.1 & 5 & \multicolumn{1}{r}{0.30} & \multicolumn{1}{r}{0.35} &
\multicolumn{1}{r}{0.34}\\
&  & 10 & \multicolumn{1}{r}{0.68} & \multicolumn{1}{r}{0.64} &
\multicolumn{1}{r}{0.70}\\
&  & 20 & \multicolumn{1}{r}{1.40} & \multicolumn{1}{r}{1.27} &
\multicolumn{1}{r}{1.40}\\\cline{2-6}
& 0.2 & 5 & \multicolumn{1}{r}{0.84} & \multicolumn{1}{r}{1.29} &
\multicolumn{1}{r}{1.27}\\
&  & 10 & \multicolumn{1}{r}{1.61} & \multicolumn{1}{r}{2.29} &
\multicolumn{1}{r}{2.33}\\
&  & 20 & \multicolumn{1}{r}{3.36} & \multicolumn{1}{r}{4.92} &
\multicolumn{1}{r}{4.79}\\\hline
\end{tabular}
\label{TabMSE}%
\end{table}%

\end{center}

The MSEs of the Optimal estimator are slightly higher than the other two for
$\varepsilon=0.1,$ and lower for $\varepsilon=0.2.$ This can be explained by
the fact that for small $\varepsilon,$ the bias is less important than the
variability, and viceversa.

\section{\medskip Appendix:\ Proofs of results\label{SecApend}}

\subsection{\bigskip Proof of Theorem \ref{Theo-Opti}\label{Apend-Opti}}

For the purposes of this proof it will be more convenient to state Hampel's
problem in its equivalent \emph{dual} form, namely, to minimize the GES under
a bound $K$ on the asymptotic variance. It is known that the solution is again
given by (\ref{esteq}), where now $m$ is a decreasing function of the bound
$K.$

For given $m$ call $\widehat{\theta}_{m}$ \ the Hampel-optimal estimator given
by (\ref{esteq}). \ In this case \ (\ref{conHamop}) takes on the form

\
\begin{equation}%
{\displaystyle\sum_{k=0}^{\infty}}
\psi_{m}^{H}(\psi_{0}(k,\widehat{\theta}_{m})-c(m,\widehat{\theta}%
_{m}))p(k,\widehat{\theta}_{m})=0. \label{cent}%
\end{equation}
When the bound $K$ tends to infinity, $m\rightarrow0$ and the GES of
$\widehat{\theta}_{m}$ tends to\ its lower bound. Then to prove the Theorem it
is enough to show that there exists $m_{0}$ such that for $m\leq m_{0},$ this
estimator \ coincides with the estimator given by (\ref{eqmed}).

We will suppose that $\psi_{0}(k,\theta)$ is \ strictly \ increasing in $k.$
The proof when it is strictly decreasing is similar. Put $k_{0}^{\ast}%
(\theta)=k_{0}(F_{\theta})$ and let%
\begin{equation}
m_{0}=\frac{1}{2}\min(\psi_{0}(k_{0}^{\ast}(\theta),\theta)-\psi_{0}%
(k_{0}^{\ast}(\theta)-1,\theta),\psi_{0}(k_{0}^{\ast}(\theta)+1,\theta
)-\psi_{0}(k_{0}^{\ast},\theta)). \label{m0}%
\end{equation}
It \ will be shown that if $m\leq m_{0}$ then\bigskip%
\begin{equation}
\psi_{0}(k_{0}^{\ast}(\theta),\theta)-m_{0}\leq c(m,\theta)\leq\psi_{0}%
(k_{0}^{\ast}(\theta),\theta)+m_{0}. \label{ineq1}%
\end{equation}

Suppose that $c(m,\theta)<\psi_{0}(k_{0}^{\ast}(\theta),\theta)-m_{0}.$ Then
we have \ $\psi_{0}(k,\theta)-c(m,\theta)>m_{0}$ \ for all $k\geq k_{0}^{\ast
}(\theta),$ and hence%
\[%
{\displaystyle\sum_{k=k_{0}^{\ast}(\theta)}^{\infty}}
\psi_{m}^{H}(\psi_{0}(k,\theta)-c(m,\theta))p(k,\theta)=m(1-F_{0}(k_{0}^{\ast
}(\theta)-1))\geq\frac{m}{2}.
\]
We also have
\[%
{\displaystyle\sum_{k=0}^{k_{0}^{\ast}-1}}
\psi_{m}^{H}(\psi_{0}(k,\theta)-c(m,\theta))\ p(k,\theta)\geq-m_{0}F_{0}%
(k_{0}^{\ast}(\theta)-1)>-\ m/2,
\]
which implies
\[%
{\displaystyle\sum_{k=0}^{\infty}}
\psi_{m}^{H}(\psi_{0}(k,\theta)-c(m,\theta))p(k,\theta)>0,
\]
contradicting ( \ref{cent} ), Similarly it can be proved that \ \ we can not
have $c(m,\theta)>\psi_{0}(k_{0}^{\ast}(\theta),\theta)-m_{0}.$

Then from \ (\ref{ineq1}) \ and (\ref{m0}) we get that $\psi_{m}^{H}(\psi
_{0}(k,\theta)-c(m,\theta))\leq-m$ \ for $k<k_{0}^{\ast}(\theta)$ and
$\psi_{m}^{H}(\psi_{0}(k,\theta)-c(m,\theta))\geq m$ for $k>k_{0}^{\ast
}(\theta).$ Then it follows from (\ref{cent}) that
\begin{equation}
-mF(k_{0}^{\ast}(\theta)-1,\theta)+m(1-F(k_{0}^{\ast}(\theta),\theta
))+(\psi_{0}(k_{0}^{\ast}(\theta),\theta)-c(m,\theta))p(k_{0}^{\ast}%
(\theta),\theta)=0, \label{cent1}%
\end{equation}
or similarly
\begin{equation}
m\frac{(1-2F(k_{0}^{\ast}(\theta),\theta))+p(k_{0}^{\ast}(\theta),\theta
)}{p(k_{0}^{\ast}(\theta),\theta)}=\psi_{m}^{H}(\psi_{0}(k_{0}^{\ast}%
(\theta),\theta)-c(m,\theta)). \label{cent2}%
\end{equation}
\ From (\ref{cent2}) \ we derive%

\begin{equation}
k_{0}^{\ast}(\theta)-0.5+\frac{(0.5-F(k_{0}^{\ast}(\theta),\theta)}%
{p(k_{0}^{\ast}(\theta),\theta)}=k_{0}^{\ast}(\theta)-1+\frac{1}{2m}\psi
_{m}^{H}(\psi_{0}(k_{0}^{\ast}(\theta),\theta)-c(m,\theta)). \label{cent3}%
\end{equation}
\ Define
\[
G_{m}(\kappa,\theta)=k-1\ +\frac{2}{2m}\psi_{m}^{H}(\psi_{0}(k,\theta
)-c(m,\theta)),
\]
$\ $\ and note that according to (\ref{cent3}) for $m\leq m_{0},$ $G_{m}$ does
not depend on $m$. Then \ \ (\ref{cent3}) \ and (\ref{umeddef})\ imply that
for $m\leq m_{0}$%
\begin{equation}
\text{\textrm{umed}}(F_{\theta})=G_{m}(k_{0}^{\ast}(\theta),\theta).
\label{cent4}%
\end{equation}

If $m\leq m_{0}$ , (\ref{esteq}) is equivalent to
\begin{equation}
-mF_{n}(k_{0}^{\ast}(\widehat{\theta}_{m})-1)+m(1-F_{n}(k_{0}^{\ast}%
(\widehat{\theta}_{m})))+(\psi_{0}(k_{0}^{\ast}(\widehat{\theta}_{m}%
),\widehat{\theta}_{m})-c(m,\widehat{\theta}_{m}))p_{n}(k_{0}^{\ast}%
(\widehat{\theta}_{m}))=0, \label{eq3}%
\end{equation}
and using the same arguments that lead to (\ref{cent4}), we can prove that
(\ref{eq3}) \ is \ equivalent to
\begin{equation}
\text{\textrm{umed}}(F_{n})=G_{m}(k_{0}^{\ast}(\widehat{\theta}_{m}%
),\widehat{\theta}_{m}). \label{eq4}%
\end{equation}

Consider the estimator $\widehat{\mathbf{\theta}}$\ defined by \textrm{umed}%
$(\widehat{\theta})=$\textrm{umed}$(F_{n}).$ By (\ref{impk0}) \ we \ have
$k_{0}(F_{n})=$ $k_{0}^{\ast}(\widehat{\theta})$ and by (\ref{cent4}) we get%
\[
\text{\textrm{umed}}(F_{n})=\text{\textrm{umed}(}F_{\widehat{\theta}}%
)=G_{m}(k_{0}^{\ast}(\widehat{\theta}),\widehat{\theta}).
\]
Then (\ref{eq4}) holds and this implies that (\ref{esteq}) holds too. This
proves the Theorem.

\subsection{Proof of Lemma\ref{Lema-Lim} \label{Apend-LemaLim}}

Let $X_{n}\sim F^{\left(  n\right)  }$ and $X\sim F,$ and call $G^{\left(
n\right)  }$ and $G$ the distributions of $X_{n}+u$ and $X+u,$ respectively,
where $u$ has a uniform distribution on $[-0.5,0.5]$ independent of $X_{n}$ or
of $X.$ Then $G^{(n)}$ and $G$ have a positive density, and $G^{(n)}%
\rightarrow_{w}G$. Since \textrm{umed}$(F^{\left(  n\right)  })=\mathrm{med}%
\left(  G^{\left(  n\right)  }\right)  $ and \textrm{umed}$(F)=\mathrm{med}%
\left(  G\right)  ,$ and the median is a weakly continuous functional, the
result is shown.

\subsection{Proof of Theorem \ref{Theo-Dist-Umed}\label{Apend-Dist-Umed}}

(a) If $F(K)>0.5,$ then then \ for large $n$ we have $k_{0}(F_{n})=K,$ and
therefore%
\[
Z_{n}=n^{1/2}\left(  \frac{0.5-F_{n}(K-1)}{p_{n}(K)}-\frac{0.5-F(K-1)}{p_{0}%
}\right)  .
\]

To derive the asymptotic distribution of $Z_{n}$ we need first to calculate
that of the vector
\[
\mathbf{d}_{n}=n^{1/2}\left[  \left(  F_{n}(K-1)-F^{1}\right)  ,\left(
\ p_{n}(K)-p_{0}\right)  \right]  ^{\prime}.
\]
The vector $d_{n}$ converges in distribution to a bivariate normal
distribution with mean $(0,0)$ and covariance matrix%
\[
\left[
\begin{array}
[c]{cc}%
a & c\\
c & b
\end{array}
\right]
\]
where
\[
a=F^{1}(1-F^{1}),~b=p_{0}(1-p_{0}),\ c=-F^{1}p_{0}.
\]

Then since for large $n$
\[
Z_{n}\backsimeq-n^{1/2}\frac{F_{n}(K-1)-F^{1}}{p_{0}}-n^{1/2}\frac
{(0.5-F^{1})(p_{n}(K)-p_{0})}{p^{2}(K)},
\]
the delta method yields that $Z_{n}\rightarrow^{D}N(0,\sigma^{2})$ where
\begin{align*}
\sigma^{2}  &  =\frac{a}{p_{0}^{2}}+\frac{b(0.5-F^{1})^{2}}{p_{0}^{4}}%
+\frac{2c(0.5-F^{1})\ }{p_{0}^{3}}\\
&  =\frac{F^{1}(1-F^{1})}{p_{0}^{2}}+\frac{(1-p_{0})(0.5-F^{1})^{2}}{p_{0}%
^{3}}-\frac{2F^{1}(0.5-F^{1})}{p_{0}^{2}},
\end{align*}
and a straightforward calculation yields (\ref{sigma-umed}).

\textbf{(b) }If $F\left(  K\right)  =0.5$ it is easy to verify that
$\mathrm{umed}(F)=K.$ Therefore for large $n$ we have%
\[
k_{0}(F_{n})=\left\{
\begin{array}
[c]{ccc}%
K & \text{if} & F_{n}(K)\geq0.5\\
K+1 & \text{if} & F_{n}(K)<0.5
\end{array}
\right.
\]

We are going to calculate $Z_{n}$ in both cases. If $F_{n}(K)\geq0.5$ we have%
\begin{align*}
\mathrm{umed}(F_{n})  &  =K-0.5+\frac{0.5-F_{n}(K-1)}{p_{n}(K)}\\
&  =K+0.5+\frac{0.5-F_{n}(K)}{p_{n}(K)},
\end{align*}
and therefore%
\begin{equation}
Z_{n}=\ \frac{n^{1/2}(0.5-F_{n}(K)}{p_{n}(K)}\leq0. \label{Zn<0}%
\end{equation}

If $F_{n}(K)<0.5$ it follows in the same way that
\[
Z_{n}=\frac{n^{1/2}(0.5-F_{n}(K)}{p_{n}(K+1)}>0.
\]

Note that, conversely, $Z_{n}\leq0$ implies $F_{n}(K)\geq0.5$ and $Z_{n}>0$
implies $F_{n}(K)<0.5.$ Since \ $n^{1/2}(0.5-F_{n}(K)\rightarrow
_{d}N(0,0.25),$ the Central Limit Theorem and Slutsky's Lemma yield
(\ref{defH}).

\subsection{Proof of Theorem \ref{Teo_Dist-Opti}\label{Apend-Dist-Opti}}

Recall that $\widehat{\theta}_{n}=g^{-1}\left(  \text{\textrm{umed}}%
(F_{n})\right)  $ and $\theta=g^{-1}\left(  \text{\textrm{umed}}(F_{\theta
})\right)  ,$ Put for brevity $K=k_{0}\left(  F_{\theta}\right)  .$

\textbf{(a) }If $F_{\theta}\left(  K\right)  >0.5$ there exists an interval
$I$ containing $\theta$ such that $t\in I$ implies that $F_{t}\left(
K\right)  >0.5$ and $F_{t}\left(  K-1\right)  <0.5,$ and therefore
$k_{0}\left(  F_{t}\right)  =K.$ Therefore $g$ is differentiable at $\theta,$
and rhe result follows from Theorem \ \ref{cons}, part (a), Theorem
\ref{Theo-Dist-Umed}, and Slutsky's Lemma.

\textbf{(b) }Assume now $F_{\theta}\left(  K\right)  =0.5.$ Then $t<\theta$
implies that $F_{t}\left(  K\right)  >0.5,$ and therefore for sufficiently
small $\delta$ we have $k_{0}\left(  F_{\theta-\delta}\right)  =K$ and
$k_{0}\left(  F_{\theta+\delta}\right)  =K+1.$ Then the left- and right side
derivatives of $g$ at $\theta$ are $g_{-}^{\prime}$ and $g_{+}^{\prime}$ given
by (\ref{laterales}), and therefore the left- and right side derivatives of
$g^{-1}$ are $1/g_{-}^{\prime}$ and $g_{+}^{\prime},$ respectively.

We have%
\begin{equation}
n^{1/2}(\widehat{\theta}_{n}-\theta)=n^{1/2}\left(  g^{-1}(\text{\textrm{umed}%
}(F_{n}))-g^{-1}(\text{\textrm{umed}}(F_{\theta}))\right)  . \label{asym1}%
\end{equation}

Note that for $\tau<0$ we have $g\left(  t\right)  -g\left(  \theta\right)
=\left(  t-\theta\right)  g_{-}^{\prime}+o\left(  t-\theta\right)  ,$ and
that
\[
\mathrm{P}\left(  n^{1/2}\left(  \text{\textrm{umed}}(F_{n}%
)-\text{\textrm{umed}}(F_{\theta})\right)  <t\right)  \rightarrow H\left(
t\right)  ,
\]
with $H$ defined in (\ref{defH}). The result follows by applying the delta method.

\medskip

\textbf{References}

Cantoni, E. and Ronchetti, E. (2001). Robust inference for generalized linear
models. \emph{Journal of the American Statistical Association}, \textbf{96}, 1022-1030.

Hampel, F. R. (1974) The influence curve and its role in robust estimation.
\emph{Journal of the American Statistical Association}, \textbf{69}, 383-394.

Ma, Y., Genton, M. and Parzen, E. (2011). Asymptotic properties of sample
quantiles of discrete distributions. \emph{Annals of the Institute of
Statistical Mathematics}, \textbf{63,} 227--243.

Maronna, R.A., Martin, R.D., Yohai, V.J. and Salibi\'{a}n-Barrera, M. (2019)
\emph{Robust Statistics: Theory and Methods (with R), Second Edition}, Wiley, Chichester.

Valdora, M. and Yohai, V.J. (1974). Robust estimators for generalized linear
models. \emph{Journal of Statistical Planning and Inference}, \textbf{146}, 31-48.

\end{document}